\documentclass[12pt]{amsart}
\usepackage{amsmath,amssymb,amsthm} 
\theoremstyle{definition}
\newtheorem{defn}{Definition}[section]
\newtheorem{remark}[defn]{Remark}
\newtheorem{example}[defn]{Example}
\theoremstyle{plain}
\newtheorem{lemma}[defn]{Lemma}
\newtheorem{proposition}[defn]{Proposition}

\newtheorem{thm}[defn]{Theorem}
\DeclareMathOperator{\im}{im}
\DeclareMathOperator{\id}{id}

\begin{document}

\title{Co-Addition for free non-associative algebras and the Hausdorff Series}


\author{Lothar Gerritzen,  Ralf Holtkamp}

\address{University of Bochum, 44780 Bochum, Germany}

\email{lothar.gerritzen@ruhr-uni-bochum.de}
\email{ralf.holtkamp@ruhr-uni-bochum.de}

\begin{abstract}
Generalizations of the series $\exp$ and $\log$ to noncommutative
non-associative and other types of algebras were considered by M.\
Lazard, and recently by V.\ Drensky and L.\ Gerritzen. There is a
unique power series $\exp(x)$ in one non-associative variable $x$
such that $\exp(x)\exp(x)=\exp(2x)$, $\exp'(0)=1$.

We call the unique series $H=H(x,y)$ in two non-associative
variables satisfying $\exp(H)=\exp(x)\exp(y)$ the non-associative
Hausdorff series, and we show that the homogeneous components
$H_n$ of $H$ are primitive elements with respect to the
co-addition for non-associative variables. We describe the space
of primitive elements for the co-addition in non-associative
variables using Taylor expansion and a projector onto the algebra
$A_0$
 of constants for the partial derivations.
 By a theorem of Kurosh,
$A_0$ is a free algebra. We describe a procedure to construct a
free algebra basis consisting of primitive elements.
\end{abstract}

\maketitle

In this article we are studying the co-addition $\Delta: K\{X\}\to
K\{X\}\otimes K\{X\}$, where $K$ is a field of characteristic 0,
and where $K\{X\}$ denotes the free unitary algebra over $K$ with
one binary operation (a non-associative, noncommutative
multiplication), also called free magma algebra generated by $X$.

The space Prim$(K\{X\})=\{f\in K\{X\}: \Delta(f)=f\otimes 1 +
1\otimes f\}$ of primitive elements for the co-addition is
considered. The first Lazard-cohomology group of the
 $\otimes$-Kurosh analyzer is given by the primitive elements for the
  co-addition, see \cite{lipseudo}, and they are also called
 pseudo-linear, following \cite{lilaloi}.
The smallest space that contains the variables and is
 closed under commutators $[f_1,f_2]:=f_1f_2-f_2f_1$ and associators
$(f_1,f_2,f_3):=(f_1f_2)f_3-f_1(f_2f_3)$ is contained in
Prim$(K\{X\})$, and is strictly smaller than Prim$(K\{X\})$. It
does not contain the primitive element $x^2x^2-2x(x^2x)+x(xx^2)$,
for example.

The algebra $(K\{X\})_0$ of constants relative to all partial
derivations $\frac{d}{dx_i}$, $x_i\in X$, contains all primitive
elements of order $\geq 2$. It is expected that $(K\{X\})_0$ is
freely generated by homogeneous primitive polynomials.

 For the
algebra $K\{x\}$ in one variable $x$, we describe a construction
of such an algebra basis for the algebra of constants, in which
Lazard-cohomology is used. The first element of the free algebra
basis consisting of primitive elements is $xx^2-x^2x$, which is
the only generator in degree 3.
 In degree
$n\geq 4$, there are $3(c_{n-1}-c_{n-2})$ generators, with Catalan
number $c_n:=\frac{(2 (n-1))!}{n!(n-1)!}$.

The basic method is based on the concept of Taylor expansion in
$K\{X\}$, which provides a projector onto the algebra $(K\{X\})_0$
of constants.

\medskip

 The classical Hausdorff series
$\log(e^xe^y)=H(x,y)=\sum_{n=1}^{\infty}H_n$, $H_n$ homogeneous of
degree $n$, is a series in associative variables with rational
coefficients.
 The famous
Campbell-Baker-Hausdorff formula states that all $H_n$ are Lie
polynomials. By a Theorem of Friedrichs, Lie polynomials are
characterized as the primitive elements of the free associative
co-addition Hopf algebra, cf.\ \cite{lireu}, p.\ 20.

Generalizations of the series $\exp$ and $\log$ to noncommutative
non-associative and other types of algebras were considered by M.\
Lazard, see \cite{lilaloi}, and recently by V.\ Drensky and L.\
Gerritzen, see \cite{lidg}. There is a unique power series
$\exp(x)$ in one non-associative variable $x$ such that
$\exp(x)\exp(x)=\exp(2x)$, $\exp'(0)=1$, and it holds that
$\exp'(x)=\exp(x)$. There exists a unique series $H=H(x,y)$
without constant term in two non-associative variables satisfying
$\exp(H)=\exp(x)\exp(y)$. We suggest to call $H$ the
non-associative Hausdorff series.

We show that the homogeneous components $H_n$ of $H$ are primitive
elements with respect to the co-addition for non-associative
variables. We obtain a recursive formula for the coefficients
$c(\tau)$ of $H$, see section 6. Each magma monomial $\tau$ is a
word with parenthesis and can be identified with an $X$-labeled
planar binary rooted tree.
 Each coefficient $\bar c(w)$,
$w$ an associative word in $x$ and $y$, of the classical Hausdorff
series is obtained as a sum  $\sum c(\tau)$ over magma monomials
$\tau$ for which the foliage is equal to $w$ (see \cite{lireu},
p.84).
 Several formulas for the coefficients in the classical case
are given in \cite{lireu} \S 3.3. There is no analogue of the
given formula, though.

\medskip

In section 1 we show that the co-addition is cocommutative,
coassociative, and that there are left and right antipodes, which
however are not anti-homomorphisms. We get a "non-associative Hopf
algebra structure" on the free magma algebra. If we identify the
monomials with $X$-labeled planar binary trees,
 the grafting of trees occurs as the multiplication, see also
\cite{licompa}, Remark (3.7).

The Hopf algebras on planar binary trees described in \cite{lilr}
or \cite{libfa} are Hopf algebra structures on free associative
algebras. The comultiplications can be described by cuts that are
given by subsets of the set of vertices.

For the co-addition Hopf algebra of section 1, the image
$\Delta(\tau)$ of a magma monomial under $\Delta$ is described in
terms of contractions of the tree $\tau$ onto subsets of the set
of leaves, see Lemma (\ref{lemdeltas}).

\medskip

Polynomials in one variable $x$ are considered in section 2. The
space of constants of degree $n$ homogeneous elements has
dimension $c_n-c_{n-1}$. We study Taylor expansions
\begin{equation*}
f=\sum_{j=0}^{\infty} x^j.a_j
\end{equation*}
of polynomials $f\in K\{x\}$ and power series $f\in K\{\{x\}\}$.
All $a_j=a_j(f)$ are uniquely defined constants for
$\frac{d}{dx}$. We define an integral
\begin{equation*}
\int f dx =\sum_{j=0}^{\infty} \frac{1}{(j+1)}x^{j+1}.a_{j}(f)
\end{equation*}
and obtain a projector
\begin{equation*}
\Phi:K\{X\}\to (K\{X\})_0
\end{equation*}
which maps $f$ onto the Taylor coefficient $a_0(f)$. In the main
result of this section, see Proposition (\ref{propmainfirst}), a
system of homogeneous free algebra generators for $(K\{X\})_0$ is
constructed. In section 4, this system is modified into a
homogeneous set of primitive polynomials, see Theorem
(\ref{thmnewbasis}).  Thus we get a co-addition Hopf algebra over
a countably infinite set of variables.

Taylor expansions for polynomials in several variables are
considered in section 3. We construct a $K$-basis for the algebra
of constants in two variables, see Proposition (\ref{propvectwo}).
In section 5 we show that the composition $f(g_1,...,g_n)$ of
pseudo-linear (i.e.\ primitive) polynomials $f, g_1,...,g_n$ is
again primitive, and that all primitive elements are obtained by
insertion of primitive elements into multi-linear ones. For
example, the primitive element $x^2x^2-2x(x^2x)+x(xx^2)$ mentioned
above can be realized as composition $f(x,x,x,x)$, where
$f(x_1,x_2,x_3,x_4)$ is the
 multi-linear primitive element
 \begin{equation*}
(x_1x_2)(x_3x_4) - x_1 (x_2(x_3x_4)) - x_4\cdot (x_1,x_2,x_3) -
x_3\cdot (x_1,x_2,x_4).
\end{equation*}
Furthermore, e.g.\
$f(x_1,[x_2,x_3],(x_4,x_5,x_6),([x_7,x_8],x_9,x_{10})]$ is also
primitive.

\medskip

We would like to thank G.\ M.\ Bergman and C.\ Mart\'inez for
pointing out references to us, and J.-L.\ Loday and V.\ Drensky
for discussions about the subject.

\begin{section}{Co-Addition.}

Let $K$ be a field of characteristic 0, and let
$X=\{x_1,x_2,...\}$ be a finite or countable set of variables. We
equip each variable $x_i$  with a non-negative degree $d_i=d(x_i)$
such that $d_{i+1}\geq d_i$.

 By $K\{X\}=F^1_{\mathcal MagAlg}(X)$ we
denote the free $K$-algebra with unit $1$ generated by $X$, called
the free magma algebra over $X$. Here a (magma) algebra over $K$
is just a vector space $V$ together with one binary operation
$\cdot: V\times V\to V$, called the multiplication (which is not
required to be associative).

 As a vector space, $K\{X\}$ has
the set of all planar binary trees with leaves labeled by the
letters $x_i$ as a basis.

The algebra $K\{X\}$ is naturally graded, $K\{X\}=\bigoplus_n
K\{X\}_n$,
  where $K\{X\}_n$ is the vector subspace generated by total degree $n$ monomials.
The set of non-trivial monomials is given by the free magma
$M=M(X)=F_{\mathcal Mag}(X)$ without unit generated by $X$, and
the total degree $d$-$\deg$ is the unique morphism $M\to {\mathbb
N}$ extending $d$.

\begin{remark}
We can embed  $K\{X\}$ into its completion
$K\{\{X\}\}=\prod_{n=0}^{\infty} K\{X\}_n$, the free complete
magma algebra generated by $X$. We have the canonical degree and order
functions $\deg$ and ord (for $d\equiv 1$).

\end{remark}

\begin{defn}\label{defcoadd}

Let $\Delta$
 be the algebra homomorphism defined on $K\{X\}$ (or $K\{\{X\}\}$)
by $x_i\mapsto x_i\otimes 1 + 1\otimes
 x_i \in K\{X\}\otimes K\{X\} \subset K\{\{X\}\}\otimes
 K\{\{X\}\}$.

The map $\Delta$ is called the co-addition.
\end{defn}

\begin{proposition}

The map $\Delta$ is coassociative,
 $(\Delta\otimes\id)\Delta=(\id\otimes\Delta)\Delta$.
It is also cocommutative, i.e.\ $\tau\circ \Delta=\Delta$,
 where $\tau$ is the permutation of tensor factors.

There is a unique $K$-linear map $\sigma$, the (left-)antipode,
such that $\sigma(w) + w +\sum \sigma(w_{<1>})w_{<2>}=0$, $\deg w
\geq 1$, where $\Delta(w)=w\otimes 1+1\otimes w +\sum
w_{<1>}w_{<2>}$.

If $\bar t$ denotes the involution induced by
$\overline{t_1t_2}=\bar t_2 \bar t_1$, then $\bar \sigma$ given by
$\bar \sigma (\bar t)=\overline{\sigma(t)}$ fulfills $\bar
\sigma(w) + w +\sum w_{<1>}\bar  \sigma(w_{<2>})=0$ ($\bar \sigma$
is the right-antipode).

\end{proposition}

\begin{proof}
To check that $(\Delta\otimes\id)\Delta=(\id\otimes\Delta)\Delta$,
we note that both homomorphisms map $x_i$ on
 $x_i\otimes 1\otimes 1 + 1\otimes
 x_i\otimes 1+ 1\otimes 1\otimes x_i$.
The antipode is recursively defined. Its properties are easy to
check.
\end{proof}

\begin{remark}

The algebra $K\{X\}$ (or $K\{\{X\}\}$) together with $\Delta$ will
be called the free (graded or complete) non-associative Hopf
algebra with co-addition. Non-associative Hopf algebras can be
defined similarly to usual Hopf algebras. (While the
comultiplication is required to be coassociative, the
multiplication is not required to be associative, and left and
right antipodes do not have to coincide). Then the co-addition is
a special comultiplication, namely the one which is dual to
addition (compare the usage in \cite{libh}, where an abelian group
valued functor leads to a co-addition).

The notion of a non-associative Hopf algebra  should not be
confused with the notion of a non-associative coalgebra in
\cite{liac}, \cite{ligr}, and \cite{lizh}, which means a not
necessarily coassociative coalgebra. The graded linear dual of
$(K\{X\}, \Delta)$ is such a coalgebra. Its non-coassociative
comultiplication is given by $t\mapsto t\otimes 1+1\otimes t +
t_1\otimes t_2$, if $t=t_1 t_2$ with $\deg(t_i)\geq 1$.

\end{remark}

\begin{remark}

Let $f=f(x_1,x_2,...)$ be a homogeneous (of total degree $\geq 1$,
say) polynomial.
 For $X'$ a copy of $X$, with
bijection $x_i\mapsto x_i'$, the substitution $x_i\mapsto
x_i+x_i'$ defines a polynomial $f(x_1+x'_1,x_2+x'_2,...)-
f(x_1,x_2,...)-f(x'_1,x'_2,...)$, which is zero iff $f$ is linear.

If we only require that $f(x_1+x'_1,x_2+x'_2,...)-
f(x_1,x_2,...)-f(x'_1,x'_2,...)=0$ holds in the case where all
variables from $X'$ commute with the variables from $X$, we get
the following weaker version of linearity.

\end{remark}

\begin{defn}\label{defpseudolin}

An element $f$ of the free magma algebra $A=K\{X\}$ (or
 its completion $K\{\{X\}\}$) is called
 ($\otimes$-)pseudo-linear, iff $f$ is a primitive element in the
free non-associative co-addition Hopf algebra, i.e.\ iff
$\Delta(f)=f\otimes 1+1\otimes f$.

The kernel of the $K$-module homomorphism given by
$\partial_1:=\Delta-\iota_1-\iota_2,$ where
$\iota_1(x_i)=x_i\otimes 1$, $\iota_2(x_i)=1\otimes x_i$, is
denoted by Prim$(A)$.

\end{defn}

\begin{remark}
\begin{itemize}
\item[(i)]
For pseudo-linear $f$, $\sigma(f)=-f$. It is worth to note that
$\sigma$ is not the anti-homomorphism induced by $x_i\mapsto
-x_i$.
\newline
 For
example, let $A=K\{x\}$. While $\sigma(x)=-x$, $\sigma(x^2)=x^2$,
we recursively compute that $\sigma(xx^2)=2xx^2-3x^2x$,
$\sigma(x^2x)=3xx^2-4x^2x$.
 The order of $\sigma$
is infinite.
\item[(ii)]
The notion pseudo-linearity follows \cite{lilaloi}. The (tuples
of) pseudo-linear elements of (\ref{defpseudolin}) form the first
cohomology group
 $H^1(\otimes-\widehat{\mbox{$\mathcal K$}})$
  of \cite{lipseudo}.
\end{itemize}

\end{remark}

\begin{defn}\label{defcontract}

Let $D_i$, $i\in {\mathbb N}$, be the unique $K$-linear derivation
defined on $A=K\{\{X\}\}$ such that $D(x_i)=1$, and $D(x_j)= 0$
else.

\smallskip

 For $s$ a tree  monomial, and $f\in K\{\{X\}\}$, we
define $\Delta_{s}(f)$ by $\Delta(f)=\sum_s \Delta_s(f)\otimes s$.

\medskip

For $t=t_1t_2$ a tree monomial, let $I(t)=I(t_1)\uplus I(t_2)$ be
the set of leaves of $t$. For any subset $I$ of $I(t)$ the
contraction $t\vert I$ is the tree with set of leaves $I$
recursively defined by $t_1t_2\vert I= (t_1\vert I\cap
I(t_1))\cdot (t_2\vert I\cap I(t_2))$, where $t\vert \emptyset=1$.
The neutral element $1$ for the grafting multiplication in $M$ is
the empty tree.

Informally the contraction $t\vert I$ is obtained by removing all
leaves not in $I$ followed by the necessary edge-contractions to
get a binary tree.

\medskip

Let $\mu_s(t)$ be the number of subsets $I$ of $I(t)$ which yield
$s$ by contraction, i.e.\ for which $t\vert I=s$.

If $t\vert I=s$, we call $s'=t\vert I^c$ a  complement of $s$ in
$t$, where $I^c$ is the complement of $I$ in $I(t)$.

\end{defn}

\begin{lemma}\label{lemdeltas}

 The following formulas hold:
\begin{itemize}
\item[(i)] For $t$ a tree, $\displaystyle
\Delta(t)=\sum_{I\subset I(t)} (t\vert I)\otimes (t\vert I^c)$;
and $\Delta_s(t)$ is given by a sum $\sum_{s'} \mu_{s'}(t)s'$ over
trees $s'$ which are complements of $s$ in $t$.
\item[(ii)]
 $\Delta_{x_i}=D_i$.
\item[(iii)]
$\Delta_s\circ\Delta_t=\Delta_t\circ\Delta_s$, for all $s,t\in M$.
\item[(iv)]
 $\Delta_{s}(fg)=
\Delta_s(f)g+f\Delta_s(g)+\Delta_{s_1}(f)\Delta_{s_2}(g)$, for
$s=s_1s_2$ of degree $>1$.
\newline
Furthermore,
$\mu_s(t)=\mu_s(t_1)+\mu_s(t_2)+\mu_{s_1}(t_1)\mu_{s_2}(t_2)$, for
$t=t_1t_2$ a tree in $M$.
\item[(v)]
Let $M_n(\{x_i\})$ be the set of elements of degree $n$ in
$M(\{x_i\})$.
\newline
 Then
$\sum_{s \in M_n(\{x_i\})}\Delta_s=\frac{1}{n!}D_i^n.$

\end{itemize}

\end{lemma}
\begin{proof}
\begin{itemize}
\item[1)] The formula in (i) follows by a direct inspection of the
co-addition map. Then (ii) and (iii) are easy consequences of (i).
Since $\Delta$ is an algebra homomorphism (and cocommutative,
coassociative), we get (iv).
\item[2)] To show (v), we proceed by induction on $\deg(t_1), \deg(t_2)$:
\newline
The sum $\ \displaystyle n!\cdot\sum_{s \in
M_n(\{x_i\})}\Delta_s(t_1t_2) =\sum_k \sum_{s_1 \in
M_k(\{x_i\}),s_2\in
M_{n-k}(\{x_i\})}\Delta_{s_1}(t_1)\Delta_{s_2}(t_2)$
\newline
 is then equal to $\sum_k {n \choose
k}D_i^kD_i^{n-k}(t_1t_2)=D_i^n(t_1t_2)$. Hence (v) follows.
\end{itemize}
\end{proof}
\bigskip

\begin{lemma}\label{critpseudo}
For $f$ homogeneous of degree $n$, $f$ is pseudo-linear if and
only if  $\Delta_s(f)=0$ for all $s\in M$ with $1\leq \deg s <
\frac{n+1}{2}$.
\end{lemma}
\begin{proof} Using cocommutativity, the criterion follows.
\end{proof}

\end{section}

\bigskip
\bigskip
\goodbreak

\begin{section}{Algebra of Constants and Taylor Expansion in one variable.}

Let $X=\{x\}$, and let $A=K\{X\}$ (or $A=K\{\{X\}\}$). Let $A_0$
be the subalgebra of elements $a\in A$ with $D(a)=0$. The algebra
$A_0$ is called algebra of constants.

By Lemma (\ref{critpseudo}), the pseudo-linear elements of order
$\geq 2$ form a subspace of $A_0$. We are going to describe $A_0$
next.

\begin{lemma}\label{lemrules}
Let $L:A\to A$ be the left multiplication by $x$, given by
$L(f):=x\cdot f$. It holds that
\begin{equation*}
\begin{split}
& D\circ L^k  - L^k\circ D = kL^{k-1}, k \geq 1\\
& D^k\circ L  - L\circ D^k = kD^{k-1}, k \geq 1.\\
\end{split}
\end{equation*}
\end{lemma}
\begin{proof}
As in the classical proof for associative variables, the equation
$D\circ L  - L\circ D = \id$ follows from the Leibniz rule. The
generalization for $k\geq 1$ follows by induction on $k$.

\end{proof}

\begin{remark}

In Proposition (\ref{proptaylor}), we are going to use the
following ordering on $M$ induced by
 the order $x_1<x_2<...$ of variables.

 Let us first order by increasing total degree. Then, if
$s=s_1s_2$ and $t=t_1t_2$ are of the same degree, we set $s<t$ if
$s_1=t_1$ and $s_2<t_2$, or if $s_1<t_1$. We will call the maximal
monomial of a homogeneous polynomial $f$ the leading term, denoted
by $f^<$.

\medskip

We write $x^j.f$ for $L^j f$, and $x^j:=L^j(1)$.

 We denote $Df$ by
$\frac{d}{dx}f$.

\end{remark}

\begin{proposition}\label{proptaylor}
\begin{itemize}
\item[(i)]
There is a unique Taylor expansion  $f=\sum_{j=0}^{\infty}
x^j.a_j$ for every $f\in A$, with $a_j=a_j(f)\in A_0$. Moreover,
the elements $a_j$ are homogeneous of degree $n-j$ if $f$ is
homogeneous of degree $n$.
\item[(ii)]
If $f=\sum_{j=0}^{\infty} x^j.a_j$, then the Taylor expansion of
$x.f$ is given by $\sum_{j=1}^{\infty} x^j.a_{j-1}$.
\item[(iii)]
The operator $\Phi$ given by $f\mapsto a_0(f)$ is a projector onto
$A_0$ with $\ker\Phi=\im L$.
\end{itemize}
\end{proposition}

\begin{proof}
\begin{itemize}
\item[1)]
We note that $D\sum_{j=0}^{\infty} x^j.a_j=\sum_{j=0}^{\infty}
(j+1)x^j.a_{j+1}$, by (\ref{lemrules}), and since $D(a_j)=0$ by
construction.
\newline
We prove the uniqueness and existence of the Taylor expansion (i)
by induction on the degree of $f$. The case $n=\deg f=0$ is
trivial. For $n\geq 1$, let $\frac{d}{dx}f$ be given by the unique
Taylor expansion $\sum_{j=0}^{\infty} x^j.b_j$. Let $g$ be given
by $\sum_{j=0}^{\infty} x^{j+1}.a_{j+1}$, where
$a_{j+1}:=\frac{b_j}{j+1}$. Then $\frac{d}{dx}(f-g)=0$, thus
$a_0:=f-g \in A_0$. Now $\sum_{j=0}^{\infty} x^j.a_j$ is the
desired Taylor expansion of $f$. For homogeneous $f$, we get
homogeneous $a_j$.
\item[2)]
Assertion (ii) follows directly from (i). For (iii), it remains to
show that $\ker\Phi\subseteq \im L$. Let us assume that $f\in
\ker\Phi$ is not in $\im L$. Then with respect to the ordering we
use, the leading term of $f\in \ker\Phi $ cannot be of the form
$x.h$. But now we easily get the contradiction $\Phi(f)^<=f^< \neq
0$.

\end{itemize}
\end{proof}

\begin{defn}

We define the integral $\int f dx$ by $\sum_{j=0}^{\infty}
\frac{1}{(j+1)}x^{j+1}.a_{j}(f)$.

\end{defn}

\begin{remark}\label{remcoeffs}
Clearly the expressions for $f$ and $\int (\frac{d}{dx}f) dx$
differ exactly by $a_0(f)$.

\smallskip

 Let us describe an algorithm to obtain the Taylor
coefficients.

For $f$ a polynomial, and $n$ maximal such that $D^nf\neq 0$,
 let $a_n=\frac{1}{n!}D^nf$. Then $\tilde f := f-x^n.a_n$ can be
 used to obtain the coefficients $a_j$, $j<n$, and we have that $D^n\tilde
 f=0$. Repeating the step, the coefficient $a_{n-i}$ can be obtained
 by applying $\frac{1}{(n-i)!}D^{n-i}$ to
 $(\id-\frac{1}{(n-i+1)!}L^{n-i+1}D^{n-i+1})\cdot...\cdot(\id-\frac{1}{n!}L^nD^n)f$.

\end{remark}

\begin{proposition}

For $f\in A$ homogeneous of degree $n$, we have

\begin{equation*}
\int f dx = \sum_{k=1}^{n+1} \frac{(-1)^{k-1}}{k!} L^k(D^{k-1} f)
\end{equation*}
This sum is finite for homogeneous $f$, and we get a continuous
operator $\int (\ ) dx:\ K\{\{X\}\}\to K\{\{X\}\}$.

\end{proposition}

\begin{proof}
Applying $D$ to the operator in question, we get
\begin{equation*}
\begin{split}
&\sum_{k=1}^{n+1} \frac{(-1)^{k-1}}{k!} D L^kD^{k-1} f=
\sum_{k=1}^{n} \frac{(-1)^{k-1}}{k!}
 L^k\circ D^k f
+\sum_{k=1}^{n+1} \frac{(-1)^{k-1}}{(k-1)!} L^{k-1}D^{k-1} f
 \\
& =\sum_{k=1}^{n+1} \frac{(-1)^{k-1}}{k!} L^k(D^{k-1} f)
-\sum_{k=0}^{n} \frac{(-1)^{k-1}}{k!} L^{k}D^{k} f = f.\\
\end{split}
\end{equation*}
Since $\sum_{k=1}^{n+1} \frac{(-1)^{k-1}}{k!} L^k(D^{k-1} f)$ is
of the form $f+ x.h$, its Taylor expansion has the same constant
term as the Taylor expansion of $\int f$, by Proposition
(\ref{proptaylor})(ii).

\end{proof}

\begin{example}

Let $f=x^2x$. Then 3 is the maximal $n$ with $D^nf\neq 0$. Thus we
set $a_3=\frac{1}{3!}D^3f=1$.

For $\tilde f:=f-x^3.a_3=f-xx^2$ we repeat the step and find that
$\tilde f\in A_0$. Thus the Taylor expansion of $f=x^2x$
 is
given by $a_0=a_0(x^2x)=x^2x-xx^2$, $a_3=1$, $a_j=0$ else.

While there are no nonzero elements in $(A_0)_2$, the space
$(A_0)_3$ is one-dimensional with generator $a_0(x^2x)$.

\end{example}

\begin{defn}

Let $\Gamma:=M-\bigl(x M\cup\{x\})$. It is a sub-magma of $M$.

\end{defn}

\begin{proposition}\label{propvecbas}
 A vector space basis of $A_{0}$
is given by the homogeneous polynomials $\Phi(s)$, $s\in \Gamma$.
The dimension of $(A_0)_n$ is given by $c_n-c_{n-1}$, where the
$c_n:=\frac{(2 (n-1))!}{n!(n-1)!}$ are the Catalan numbers.

\end{proposition}

\begin{proof}

Let $M_n$ be the subset of $M$ consisting of degree $n$ monomials,
and let $\Gamma_n$ be the subset of $\Gamma$ of elements of degree
$n$. Then the number of elements in $\Gamma_n$ is given by
$\#\Gamma_n=\# M-\# (x M\cup\{x\}) =c_n-c_{n-1}$.

By Proposition (\ref{proptaylor})(iii), we get a vector space
basis $\Phi(s)$, $s\in \Gamma$.

\end{proof}

\bigskip

\begin{proposition}\label{propomega}
\begin{itemize}
\item[(i)]
The magma $\Gamma$ has the following set as a free generating set:
\begin{equation*}
\begin{split}
 \Omega &=\{ s_1\cdot s_2 :\ \ s_1\in M-\{x\},s_2\in M, \text{
s.t.\ } s_1\not\in \Gamma \text{ or } s_2\not\in \Gamma \}
\\
&=\{ (x^i.v_1)\cdot (x^j.v_2) :\ \ v_1,v_2\in \Gamma\cup 1,
i+j\geq 1, v_1\neq 1 \text{ if } i\leq 1, v_2 \neq 1  \text{ if }
j=0 \}\\
\end{split}
\end{equation*}
\item[(ii)]
Let $\Omega_n:=\{\omega\in \Omega: \deg(\omega)=n\}$. We have that
$\Omega_1=\Omega_2=\emptyset$, $\Omega_3=\{x^2x\}$.
\newline
For $n\geq 4$, $\# \Omega_n=3(c_{n-1}-c_{n-2})$.
\end{itemize}
\end{proposition}

\begin{proof}
\begin{itemize}
\item[1)]
Clearly, $\Gamma-\Gamma\cdot \Gamma$ is a generating set. Since
the restriction of the multiplication $M\times M\to M$ on $\Gamma$
is injective as well, the given generating set is free.
\item[2)]
The set $\Omega_n$ is the union $\Omega'_n\cup \Omega''_n$ of
$\Omega'_n=\{(x \cdot t_1) t_2: t_1,t_2\in M,
\deg(t_1)+\deg(t_2)=n-1\}$ and $\Omega''_n=\{t_1 (x\cdot t_2),
t_1\in M-\{x\}, t_2\in M\cup\{1\},\deg(t_1t_2)=n-1\}$.
\newline
The number of elements of $\Omega'_n$ is given by $c_{n-1}$, since
the pairs $(t_1,t_2)$ can be identified with trees $t_1t_2$ of
degree $n-1$. For $\Omega''_n$, one similarly counts $c_{n-1}$
elements of the form $t_1\cdot x$ plus $c_{n-1}-c_{n-2}$ elements
of the form $t_1(x\cdot t_2)$ with $t_1,t_2 \in M, t_1\neq x$.
\newline
The number of elements of $\Omega'_n\cap
\Omega''_n=\{(xt_1)(xt_2): t_1,t_2\in M,\deg(t_1t_2)=n-2 \}\cup
\{(xt_1)x: t_1\in M , \deg(t_1)=n-2 \}$ is given by $2c_{n-2}$.
Thus assertion (ii) follows.
\end{itemize}
\end{proof}

\begin{proposition}\label{propmainfirst}
Let $E=\{\Phi(\omega):\omega\in \Omega\}$. Then $E$ is a sequence
\begin{equation*}
y_{3,1}=a_0(x^2x),\  y_{4,1},...,y_{4,3c_3-3c_2},\
y_{5,1},...,y_{5,3c_4-3c_3},\ ...
\end{equation*}
 of elements in
$A_0$, ordered by the leading monomials as in (\ref{proptaylor}),
that freely generates the algebra $A_0$.
\end{proposition}

\begin{proof}
 Since $\Phi$ does not change the leading term, the assertion
 follows from Proposition  (\ref{proptaylor})(iii)
 and Proposition (\ref{propomega}).
\end{proof}

\begin{example}\label{extaylorfirst}

The tree $(x^2 x)^2\in \Gamma$ is not an element of $\Omega$. It
is the smallest element of $\Gamma$ that is not contained in
$\Omega$.

To determine the Taylor expansions of the  trees $t=t_1t_2$ of
degree $\geq 4$ in $\Omega$, let us recall that by
(\ref{lemdeltas})(iv) and (v),
$\mu_s(t)=\mu_s(t_1)+\mu_s(t_2)+\mu_{s_1}(t_1)\mu_{s_2}(t_2)$ is
the coefficient of $s$ in $\frac{1}{n!}D^n$.

For $s=x^2x$, the formula reads
$\mu_{x^2x}(t)=\mu_{x^2x}(t_1)+\mu_{x^2x}(t_2)+{\deg t_1 \choose
2}\deg t_2$.

\smallskip

While $a_n=1, a_{n-1}=a_{n-2}=0,$ the Taylor coefficient
$a_{n-3}=\frac{1}{(n-3)!}D^{n-3}(f-x^{n})$ can be determined by
this formula, because it has to be of the form $\alpha\cdot
y_{3,1}$ with $\alpha\in{\mathbb Q}$. It is simply given by
$\mu_{x^2x}(t)$, as $(D^{n-3} x^n)^<=x^3$.

In degree $4$, the Taylor expansions are given by coefficients
$a_4=1, a_3=a_2=0,$ and $a_1, a_0$ as follows:

For $x^2x^2$, $a_1=2y_{3,1}$, $a_0=x^2x^2-2x.y_1-x^4.1
=
x^2x^2-2x(x^2x)+x^4$.

For $x^3x$, $a_1=3y_{3,1}$, $a_0=x^3x-3x.y_1-x^4
=
x^3x-3x(x^2x)+2x^4$.

For $((x^2)x)x$, $a_1=4y_{3,1}$, $a_0=((x^2)x)x-4x.y_1-x^4
=((x^2)x)x-4x(x^2x)+3x^4$.

The elements $y_{4,1}=a_0(x^2x^2), y_{4,2}=a_0(x^3x),
y_{4,3}=a_0(((x^2)x)x)$ form a basis of $(A_0)_4$.

\end{example}

\bigskip
\begin{remark}

Similarly, for $t$ a tree of degree $n$, and $1\leq r\leq n-3$,
the Taylor coefficient $a_r(t)$ is of the form
$\sum\alpha_{n-r,i}\cdot y_{n-r,i}$, with
$\alpha_{n-r,i}=\mu_{y_{n-r,i}^{<}}(t) \in {\mathbb N}$.

\end{remark}

\begin{example}

For $s=x^lx^{r-l}$ and $t=x^kx^{n-k}$ elements of $\Omega$ with
$r=\deg s\geq 3$, $n=\deg t\geq3$, and such that $2\leq l \leq
r-1$, $2\leq k \leq n-1$, one easily shows that $\mu_s(t)={k
\choose l}{n-k \choose r-l}$.

Thus $z_{n,k}:=a_0(x^kx^{n-k})$ is given by
$x^kx^{n-k}-x^n-\sum_{z_{r,l}}{k \choose l}{n-k \choose r-l}
x^{n-r}.z_{r,l}$, where the sum ranges over all elements $z_{r,l}$
with leading terms $x^lx^{r-l}$.

E.g., for $x^2x^3$, we get $a_2=3y_1$, $a_1=3y_2$,
$z_{5,2}=x^2x^3-x^5 -3x^2.z_{3,2}- 3x.z_{4,2}$.

\smallskip

In degree 5, $(A_0)_5$ is 9-dimensional. To obtain a basis we take
$z_{5,2}, z_{5,3}, z_{5,4}$, and we take the $a_0$-terms of the
following six trees, the Taylor coefficients of which we give
here: The tree $x^2(x^2x)$, with coefficients $a_2=4y_{3,1}$,
$a_1=3y_{4,1}$; the tree $(x^2x)x^2$, with $a_2=7y_{3,1}$,
$a_1=3y_{4,1}+2y_{4,3}$; the tree $(x (x^2x))x$, with
$a_2=7y_{3,1}$, $a_1=3y_{4,2}+y_{4,3}$; the tree $(x^2x^2)x$, with
$a_2=8y_{3,1}$, $a_1=y_{4,1}+2y_{4,2}+2y_{4,3}$; the tree
$(x^3x)x$, with $a_2=y_{3,1}$, $a_1=2y_{4,2}$; and the tree
$((x^2x)x)x$, with $a_2=10y_{3,1}$, $a_1=5y_{4,3}$.

\end{example}

\end{section}

\bigskip
\bigskip

\begin{section}{Taylor Expansion for Several Variables.}

Let $X=\{x_1,x_2,...\}$ and $A=K\{X\}$ (or $A=K\{\{X\}\}$). Let
$L_i$ denote the left multiplication by $x_i$. For ${\bf j}$ a
tuple $(j_1,...,j_r)$ with entries in ${\mathbb N}$, we set
$x^{\bf j}.f:= \bigl(L_1^{j_1}\circ\ ... \circ
L_r^{j_r}\bigr)(f)$.

\begin{proposition}\label{proptayplus}
\begin{itemize}
\item[(i)]
Let $B$ be a graded magma algebra with unit, and let $x\in B$. Let
a derivation $D:B\to B$ be given, such that $D$ is a nilpotent
operator with $D(x)=1$. Then there is a unique Taylor expansion
$f=\sum_{j=0}^{\infty} x^j.b_j$ for every $f\in B$, with
$D(b_j)=0$ for all $j$.
\item[(ii)]
For each $x_i\in X$, there is a unique Taylor expansion
$f=\sum_{j=0}^{\infty} x_i^j.b_j$ (with respect to one variable)
for every $f\in A$, with $\frac{d}{dx_i}b_j=0$ for all $j$.
\newline
If $\Gamma(X,i)$ denotes the sub-magma of $M(X)$ given by
$M(X)-(x_i M(X)\cup \{x_i\})$, then the elements $b_0(s), s\in
\Gamma(X,i)$, form a vector space basis of the algebra of
constants with respect to $D_i$.
\item[(iii)]
There is a unique (total) Taylor expansion  $f=\sum_{\bf j} x^{\bf
j}.a_{\bf j}$ for every $f\in  A$, such that all $a_{\bf j}$ are
in $A_0:= \{f\in A: \frac{d}{dx_i}f=0$, all $i \}$.
\end{itemize}
\end{proposition}

\begin{proof}
Assertion (i) is proven as in Proposition (\ref{proptaylor}),
where the degree of $f$ is now replaced by the smallest $n$ such
that $D^n f=0$.
\newline
Using (i) we can prove assertion (ii) as in the case of one
variable.
\newline
Having obtained the Taylor expansion with respect to the variable
$x_1$, we can expand all coefficients again with respect to the
next variable. The set $\tilde
X:=\Omega(\Gamma(X,1)):=\Gamma(X,1)-\Gamma(X,1)\cdot \Gamma(X,1)$
is a free generating set for the algebra of constants of
$\frac{d}{dx_1}$. Iterating the process, we get as the result an
expansion $f=
 \sum_{\bf j} x_1^{j_1}.(x_2^{j_2}.(...\ .a_{\bf j})...)=
 \sum_{\bf j} x^{\bf j}.a_{\bf j}$.
The elements $a_{\bf j}$ are constants with respect to all
$\frac{d}{dx_i}$.

\end{proof}

\begin{remark}
If $\Phi_i$ is given by $\Phi_i(f)=f-\int D_i f dx_i$, then
\begin{equation*}
\Phi_i = \sum_{k=0}^{\infty} \frac{(-1)^{k}}{k!} L_{i}^k\circ
D_i^{k}.
\end{equation*}

Let $\Phi: f\mapsto a_{(0,...,0)}$ be the projector $A\to A_0$
obtained by the above Taylor expansion of several variables.

\end{remark}

\begin{proposition}\label{propvectwo}

For $X=\{x,y\}$, a vector space basis of
$A_0=\im\Phi=\im(\Phi_2\circ\Phi_1)$ is given by all $\Phi(s)$,
$s\in \Gamma$, where
\begin{equation*}
\Gamma=\bigl\{s_1s_2\in M(x,y): \deg s_1\geq 2,\deg s_2\geq
1\bigr\}\cup\bigl\{y\cdot(x \cdot t) \in M(x,y): \deg t\geq
0\bigr\}.
\end{equation*}

\end{proposition}

\begin{proof}
\begin{itemize}
\item[1)]
First we extend the basis of $\ker\Phi_1=\{x\cdot h: h\in A\}$,
see Proposition (\ref{propvecbas}),
 to a basis
of $\ker(\Phi_2\circ\Phi_1)$. To obtain the desired basis
elements, we can without loss of generality look at the space of
all $f\in A$ with $f=\Phi_1(f)=y\cdot g$ for some $g\in A$. Then
necessarily $\frac{d}{dx}(g)=0$. A basis for this space is the set
$\{y\Phi_1(t): x\neq t\in M-xM\}$.
\item[2)]
 To get a basis for $\im\Phi$, we can include all $\Phi(s)$ with
$s=s_1s_2\in M(x,y)$ with $\deg s_1\geq 2$. When we furthermore
take elements $\Phi_2(\Phi_1(s))$ we have to exclude the trees $s$
with $\Phi_1(s)^<=y\Phi_1(t)^<$, by 1). For $t\not\in x(M(x,y)\cup
\{1\})$, $y \Phi_1(t)^<=\Phi_1(y t)^<$, and the assertion follows.

\end{itemize}
\end{proof}

\begin{example}\label{extwotayl}

Let $X=\{x,y\}$. The monomials $xy$ and $y(yx)$ are not in
$\Gamma$.

 The monomial $t=yx$ is an element of $\Gamma$. Its
Taylor expansion $\sum_{j=0}^{\infty} x^j.b_j$ with respect to $x$
is given by $x.b_1+b_0$ with $b_1=y$, $b_0=yx-xy$. Applying Taylor
expansion with respect to $y$ on $b_1$, $b_0$, we get the (total)
Taylor expansion $t= x^{(1,1)}.1 + x^{(0,0)}.a_{(0,0)}(t)$ with
$a_{(0,0)}(yx)=yx-xy.$

The Taylor expansion of $y(xy)\in \Gamma$ with respect to $x$ is
given by $b_1=y^2, b_0=y(xy)-xy^2$. Its total Taylor expansion is
$t= x^{(1,2)}.1 + x^{(0,1)}.a_{(0,0)}(yx)+ a_{(0,0)}(t) $ with
$a_{(0,0)}(t)=2y(xy)-xy^2-y(yx)$. We note that the coefficient of
$y(xy)$ is not 1 in $a_{(0,0)}(y(xy))$, and that $y(xy)$ is not
the leading monomial.

For $t=y(y(yx))\in \Gamma$, we get
$a_{(0,0)}(t)=y(y(yx))-3y(y(xy))+3y(xy^2)-x(yy^2)$.

\end{example}

\begin{defn}\label{defxupper}

For $X=\{x_1,x_2,...\}$ a set of variables, let

 $x_j^{(p)}\in
\underbrace{K\{X\}\otimes ... \otimes K\{X\}}_{\geq p}$ denote
$\underbrace{1\otimes ... \otimes 1}_{p-1}\otimes x_j \otimes 1
\otimes ... \otimes 1$.

We regard the tensors $x_j^{(p)}$ as variables (and note that some
associativity and commutativity occurs among them).

\end{defn}

Let $X=\{x\}$,  let $n\geq 2$, and let ${\bf j}\in {\mathbb N}^n$.

We will need the following variation of Proposition
(\ref{proptayplus}) in order to express the formula in Proposition
(\ref{propaoo}). Here $x^{\bf j}.f$ denotes the image of $f\in
A^{\otimes n}$ under $L^{j_1}\otimes...\otimes L^{j_n}:A^{\otimes
n}\to A^{\otimes n}$.

\begin{proposition}\label{proptwotayl}
\begin{itemize}
\item[(i)]
There is a unique Taylor expansion $f=\sum_{\bf j} x^{\bf
j}.a_{\bf j}$ for every $f\in A^{\otimes n}$, with $a_{\bf j}\in
A_0^{\otimes n}$.
\item[(ii)]
The operator $A^{\otimes n}\to A_0^{\otimes n}, f\mapsto
a_{0,...,0}$, is a projector with $\ker\Phi=\bigcup\im L_i$ and
equal to $\Phi^{\otimes n}$, where $\Phi:K\{x\}\to K\{x\}_0,
f\mapsto a_0(f)$.
\end{itemize}
\end{proposition}

\begin{proof}
The Taylor expansion with respect to variables $x^{(p)}, 1\leq p
\leq n$, is defined completely analogously to Proposition
(\ref{proptayplus}).
\newline
Having expanded with respect to $x^{(p)}$, the resulting expansion
 has coefficients in
\newline
  $A\otimes\ ...\otimes A_0\otimes A \otimes
... \otimes A$ ($A_0$ in $p$-th position).
\newline
We can expand the coefficients again by the remaining $x^{(p)}$,
and we finally get a Taylor expansion $f=\sum_{\bf j} x^{\bf
j}.a_{\bf j}$. The difference to the situation in Proposition
(\ref{proptayplus}) is that the expansion is
 independent on
the chosen order of steps. The coeffients $a_{\bf j}$ are
constants with respect to all $\frac{d}{d x^{(p)}}$, and the
projector given by $f\mapsto a_{0,...,0}$ is equal to
$\Phi^{\otimes n}$.
\end{proof}

\end{section}

\bigskip
\bigskip

\begin{section}{Primitive elements in the case of one variable.}

Let $X=\{x\}$ and $A_0$ be the algebra of constants.

 We have
constructed the generators of $\bigoplus_{i=1}^5 (A_0)_i$ in a way
that they are all pseudo-linear. The elements $\Phi(t), t\in
\Omega_{\geq 6},$ are no longer pseudo-linear in general, though.
We are going to show, that $A_0$ is a (cocommutative) Hopf
sub-algebra of the co-addition Hopf algebra.

\begin{proposition}\label{propaoo}

The map $\Delta$ restricts to an algebra homomorphism $A_0\to
A_0\otimes A_0$.

Furthermore, for $t$ a tree monomial,
\begin{equation*}
\Delta(a_0(t))=a_{0,0}(\Delta(t))=a_0(t)\otimes 1 + 1\otimes
a_0(t) + \sum_{\emptyset \neq I\subsetneq I(t)} a_0(t\vert
I)\otimes a_0(t\vert I^c).
\end{equation*}

\end{proposition}
\begin{proof}

\begin{itemize}
\item[1)]
By Lemma (\ref{lemdeltas}), (ii),(iii), we get that
\begin{equation*}
 \Delta\circ D  = (D\otimes \id)\circ \Delta= (\id\otimes D)\circ
\Delta.
\end{equation*}
Hence if $Df=0$ then $\Delta(f)\in (\ker D\otimes A) \cap (A
\otimes \ker D)$.
\item[2)]
For $f=\sum x^k.a_k$, $\Delta(f)=\sum\sum {k \choose i}
(x^i\otimes x^{k-i}).\Delta(a_k)$.
\newline
By Proposition (\ref{proptwotayl}),
$a_{0,0}(\Delta(f))=\Phi^{\otimes
2}(\Delta(f))=\Phi(\Delta(a_0))$. Now
$\Phi(\Delta(a_0))=\Delta(a_0),$ as $\Delta(a_0)\in A_0\otimes
A_0$.
\end{itemize}
\end{proof}

\begin{remark}

Lazard-Lie theory (see \cite{lilaloi}), in the generalized setting
presented in \cite{lipseudo} (involving the tensor product of not
necessarily commutative or associative algebras), can be used to
show that $A_0$ is again the free non-associative Hopf algebra
with co-addition (over an infinite set of generators).

We fix some notation first.

\end{remark}

\begin{defn}

Let $X=\{x_1,x_2,...\}$ be a set of variables, and let
 $x_j^{(p)}$ be defined as in Definition (\ref{defxupper}).

Let $\bigl\vert\otimes\text{-}{\mathcal MagAlg}(X)\bigr\vert^{n}$
denote the subalgebra without unit generated by all variables
$x_j^{(p)}, 1\leq p \leq n$, in $K\{X\}^{\otimes n}$. For $1\leq
p\leq n$, let $\partial_p$ be the
 $K$-linear map
$\underbrace{\id\otimes\ ...\otimes\id}_{p-1}\otimes\
\partial_1 \otimes\id \otimes\ ...\otimes
\id: \bigl\vert\otimes\text{-}{\mathcal MagAlg}(X)\bigr\vert^{n}
\to\bigl\vert\otimes\text{-}{\mathcal MagAlg}(X)\bigr\vert^{n+1}$.

\smallskip

Generalizing the $n=1$ case, given $1\leq p \leq n$, an element
$f$ of $\bigl\vert\otimes\text{-}{\mathcal
MagAlg}(X)\bigr\vert^{n}$
 is called pseudo-linear with
respect to the $p$-th tensor argument, iff $\partial_p(f)=0$.

\end{defn}

\begin{thm}\label{thmnewbasis}
Let $K$ be a field of characteristic 0. The algebra $A_0$ of
constants is a free algebra generated by a set of variables
$Y=\{y_1,y_2,...\}$ with non-negative degrees $d_i$ such that
$d_{i+1}\geq d_i$, $d_1=3$. Let $\Delta': K\{Y\}\to K\{Y\}\otimes
K\{Y\}$ be given by $\Delta\vert A_0$.

Then $(K\{Y\},\Delta')$ is isomorphic to the free non-associative
 co-addition Hopf algebra $(K\{Y\},\Delta)$.
  The isomorphism $\varphi$ with
$\Delta=(\varphi\otimes \varphi)\circ\Delta'\circ\varphi^{-1}$ is
strict in the sense that its linear part is the identity on the
vector space generated by all $y_i$.

\end{thm}

\begin{proof}(Sketch of Proof.)
\begin{itemize}
\item[1)]
We have verified the result that $A_0$ is a free algebra, which
follows also from a theorem of Kurosh, cf.\ \cite{liku}. The
restriction $\Delta':=\Delta\vert A_0$ is a coassociative and
cocommutative algebra homomorphism, graded in the sense that the
total $d$-degree of each image $\Delta'(y_i)$ is again $d_i$.
 The co-addition from
Definition (\ref{defcoadd}) is given by the (infinite) tuple
$F=(y_1^{(1)}+y_1^{(2)},y_2^{(1)}+y_2^{(2)},...)$ of images of the
free algebra generators $y_i$.
  The augmentation map is the
counit both for $\Delta$ and $\Delta'$. Thus $\Delta'$ is given by
a tuple $F'$ whose homogeneous degree $r=1$ part (w.r.t.\ the
canoncial degree function for free variables $y_i$) is also $F$.
\newline
Restricted to the variables $y_i$ with $d$-degree $d_i\leq 5$, the
map $\Delta'$ is already the co-addition.
\item[2)]
Let now $r\geq 2$, and let us restrict to the subalgebra generated by the
set $Y^{<}=\{y_i: d_i\leq r\cdot d_1\}$. We suppose that up to
degree $r-1$ terms, we have constructed an isomorphism $\varphi$
such that $(\varphi\otimes \varphi)\circ\Delta'\circ\varphi^{-1}$
is
 the co-addition (in other words, the two associated $(r-1)$-chunks for $F$, $F'$ coincide).
\newline
The homogeneous part of degree $r$ of $F'$ (or equivalently, of
$F'-F$) is a tuple of elements from
$\bigl\vert\otimes\text{-}{\mathcal MagAlg}(Y^<)\bigr\vert^{2}$.
By Theorem (3.5) of \cite{lipseudo}, the tuples of elements of
$\bigl\vert\otimes\text{-}{\mathcal MagAlg}(Y^<)\bigr\vert^{n}$,
$n\in {\mathbb N}$, form a pseudo-analyzer.
 Since both $F$ and $F'$ yield
coassociative extensions of the given $(r-1)$-chunk,
 Proposition (4.7)(ii) of \cite{lipseudo} can be
applied. Thus the homogeneous part of degree $r$ of $F'$ is a
cocycle $E\in\ker\delta_2$, where the coboundary map $\delta_2$ is
given by $\partial_2-\partial_1$ (applied component-wise). The
isomorphism class of the extension does not depend on
coboundaries, i.e.\ images of $\delta_1:=-\partial_1$. Now $r!E$
is $0$ in $\ker \delta_2/ \im \delta_1$ (as $E$ is symmetric
w.r.t.\ all permutations). Since char$(K)=0$, there is a
homogeneous element $\alpha_r$ of degree $r$ with
$E=\delta_1(\alpha_r)$. Now $-\delta_1(\alpha_r)$ is the degree
$r$-part of the modification of $F'$ by
 the bijective transformation given by the tuple $\varphi'=\varphi +
\alpha_r$ of elements of $K\{Y\}$. Thus $(\varphi'\otimes
\varphi')\circ\Delta'\circ{\varphi'}^{-1}$ is the co-addition up
to degree $r$ terms.
\newline
Including new generators $y_i$ in $Y^{<}$ next, observe that the
necessary $(\leq r)$-modifications can be made without changing
the already modified generators (and the $\leq r$ components of their images).
\end{itemize}

\end{proof}

\bigskip

\begin{proposition}\label{propsubhopf}
The algebra freely generated by the elements
$z_{n,k}=a_0(x^kx^{n-k})$, $n\geq 3$, $2\leq k\leq n-1$, is a Hopf
subalgebra of $A_0$. The image $\Delta(z_{n,k})$ is given by
\begin{equation*}
z_{n,k}\otimes 1 + 1\otimes z_{n,k} + \sum_{m=2}^{k-2}\ \
\sum_{m+1\leq\  l\ \leq m+n-k-1}
 {k \choose m}{n-k \choose l-m}
  z_{l,m}\otimes
 z_{n-l,k-m}.
\end{equation*}
Especially, $z_{n,n-1}=a_0(x^{n-1}\cdot x)$,
$z_{n,2}=a_0(x^{2}\cdot x^{n-2}),$ and $z_{n,3}=a_0(x^{3}\cdot
x^{n-3})$ are always
 pseudo-linear.

\end{proposition}

\begin{proof}

By Proposition (\ref{propaoo}),
$\Delta(z_{n,k})=\Delta(a_0(x^kx^{n-k}))$ is given by
 $ \sum_{I}
a_0(x^kx^{n-k}\vert I)\otimes a_0(x^kx^{n-k}\vert I^c)$, which is
a sum $\sum
\mu_{x^{l_1}x^{l_2}}(x^{k}x^{n-k})z_{l_1+l_2,l_1}\otimes
z_{n-(l_1+l_2),k-l_1}$ over all possible $l_1\leq k,l_2\leq
n-l_1$.

By induction, one verifies that $\mu_{x^n}x^m={m \choose n}$ and
that

$\mu_{x^{l_1}x^{l_2}}(x^{k_1}x^{k_2})={k_1\choose l_1}{k_2 \choose
l_2}$, if $l_1,k_1\geq 2,$ and $l_2,k_2\geq 1$.

\end{proof}

\begin{example}

Let us illustrate the modification process of theorem
(\ref{thmnewbasis}) in the first nontrivial case, i.e.\ we look at
the generators $y_i$ with $d_i=6$.

By Proposition (\ref{propsubhopf}), the elements $z_{6,5},
z_{6,3}$ and $z_{6,2}$ of $(A_0)_6$ are pseudo-linear.  It is also
easy to see that $\Phi(x^3(x^2x))$ and $\Phi((x^2x)x^3)$ (given by the remaining
trees in $\Omega$ that are a product of two $\deg 3$-monomials)
are pseudo-linear.

For $z_{6,4}=\Phi(x^4x^2)$,
\newline
$\Delta(z_{6,4})=z_{6,4}\otimes 1+1\otimes z_{6,4} + {4 \choose
2}{2 \choose 1}
  z_{3,2}\otimes
 z_{3,2}=z_{6,4}\otimes 1+1\otimes
z_{6,4} +12a_0(x^2x)\otimes a_0(x^2x)$.

Now $z_{6,4}-6a_0(x^2x)\cdot a_0(x^2x)$ is pseudo-linear, because
$\Delta((a_0(x^2x))^2)=(\Delta(a_0))^2=a_0\otimes 1+1\otimes
a_0+2a_0\otimes a_0$.

For $t= (((x^2x)x)x)x$, the element  $\Phi(t)$ is not pseudo-linear, as
$\Delta(a_0(t))=t\otimes 1 +1 \otimes t + {6 \choose 3}
a_0(x^2x)\otimes a_0(x^2x)= t\otimes 1 +1 \otimes t +
20a_0(x^2x)\otimes a_0(x^2x)$. But $\Phi(t)-10(x^2x-x^3)^2$ is
pseudo-linear.

Similarly we can handle
 the remaining four trees $t_1\cdot x^2$, $\deg
t_1=4$, the twelve trees of the form $t_1\cdot x$, $t_1\neq x^5$,
and the four trees $x^2\cdot t_2$, $t_2\neq x^4$.

\end{example}

\end{section}

\bigskip
\bigskip
\goodbreak

\begin{section}{Primitive elements for Several Variables.}

Let the set $X$ of variables have $n\geq 1$ or countably many
elements.

\begin{lemma}\label{lemmaitcom}
The following elements (and their $K$-linear combinations) are
pseudo-linear (primitive for the co-addition):
\begin{itemize}
\item[(i)] variables $x_i\in X$
\item[(ii)] the commutators $[f_1,f_2]:=f_1f_2-f_2f_1$ of pseudo-linear $f_i, i=1,2$
\item[(iii)] the associators $(f_1,f_2,f_3):=(f_1f_2)f_3-f_1(f_2f_3)$ of pseudo-linear $f_i, i=1,2,3.$
\end{itemize}
\end{lemma}
\begin{proof}
Since the elements $f_i$ are pseudo-linear, it is clear that
computing the deviation to the co-additive part we can assume that
for each $i$, the leaves $x$ of $f_i$ either all have to be
substituted by the corresponding $x^{(2)}$ or all by $x^{(1)}$.
The associativity and commutativity relations for the variables
$x_j^{(1)}, x_j^{(2)}$ force the resulting expressions to be zero.
\end{proof}

\begin{remark}
 For $X=\{x\}$, the space $(A_0)_4$ is
3-dimensional, and $[x,(x,x,x)]=a_0(x^3x)-a_0((x^2x)x)$ is the
only pseudo-linear element which can be generated by the process
of (\ref{lemmaitcom}).
\end{remark}

\begin{lemma}\label{lemcompo}
\begin{itemize}
\item[(i)] For$ \# X \geq n$, let $f$ be a pseudo-linear element
in $n$ variables. Then for all pseudo-linear $g_1,...,g_n$, the
composition $\eta_{g_1,...,g_n}(f)=f(g_1,...,g_n)$ is again
pseudo-linear.
\item[(ii)]
If $\# X \geq 2$, each Lie polynomial $f=\sum x_{i_1}\cdot ...
\cdot x_{i_r}$ over $X$ in associative variables induces a
pseudo-linear element $f=\sum x_{i_1} (x_{i_2}\cdot (...
x_{i_r})...)$ with right normed brackets inserted.
\end{itemize}
\end{lemma}
\begin{proof} The algebra homomorphism $K\{X\}\to K\{X'\}$ given on the generators by $x_i\mapsto g_i$
is a homomorphism for the co-addition, if the elements $g_i$ are
pseudo-linear. Thus (i) follows. Assertion (ii) is an easy
observation.
\end{proof}

\begin{example}\label{extwovar}
The Lie polynomial
 $[y,[y,x]]=yyx+xy^2-2yxy$ in associative variables $x$ and $y$ leads to the
 pseudo-linear element $y(yx)+xy^2-2y(xy)$, which is equal to $=-a_{0,0}(y(xy))$ by Example
 (\ref{extwotayl}).
Similarly, $[x,[x,y]]$ leads to $yx^2+ x(xy)-2x(yx)$, which is
equal to $a_{0,0}(yx^2)$.

Using Proposition (\ref{propvectwo}), for $X=\{x,y\}$ a basis of
the space of pseudo-linear elements in degree 3 can be given by
\begin{equation*}
\begin{split}
&a_{0,0}(xx^2),a_{0,0}(yy^2), a_{0,0}(yx^2), a_{0,0}(y(xy)),\\
 &
a_{0,0}((xy)y)= (xy)y-xy^2=(x,y,y),\\
 & a_{0,0}((yx)y)=
(yx)y+y(xy)-y(yx)-xy^2=(y,x,y)+a_{0,0}(y(xy)),\\
 & a_{0,0}(y^2x)=
y^2x-2y(yx)+2y(xy)-xy^2=(y,y,x)+a_{0,0}(y(xy)),\\
 & a_{0,0}((yx)x)=(yx)x-x(xy)=(y,x,x),\\
 & a_{0,0}((xy)x)=(xy)x-x(xy)=(x,y,x),\\
 & a_{0,0}(x^2y)= x^2y-x(xy)=(x,x,y).
\end{split}
\end{equation*}

\end{example}

\begin{example}
 The multi-linear element $f(x_1,x_2,x_3,x_4)$ given by
\begin{equation*}
 a_{\bf 0}\bigl(
(x_1x_2)(x_3x_4)\bigr)= (x_1x_2)(x_3x_4)  - x_4\cdot (x_1,x_2,x_3)
- x_3\cdot (x_1,x_2,x_4) - x_1\cdot (x_2(x_3x_4))
\end{equation*}
 is pseudo-linear. Using (\ref{lemcompo})(i), we see that $f(x,x,x,x)$ is
pseudo-linear, too. It is equal to $a_0(x^2x^2)$, see
(\ref{extaylorfirst}), and is linearly independent of
$[x,(x,x,x)]$.

\end{example}

\begin{remark}\label{remmulti}
The co-addition $\Delta$ respects the multi-degree.

Let a homogeneous pseudo-linear element $f$ be given which is not
multi-linear. For pseudo-linear $g_1,..., g_n$, e.g.\ sums of
variables, not only $f(g_1,...,g_n)$ but also all its multi-degree
components are pseudo-linear. In this way we can obtain a
pseudo-linear multi-linear element $\tilde f$ that yields $f$ (up
to a constant factor) by evaluation on the original set of
variables, cf.\ \cite{lilaloi}, \S 2. Hence all pseudo-linear
elements can be obtained by evaluation of multi-linear ones.

Thus the construction of multi-linear primitive elements is
crucial, and we will study this problem in a subsequent paper.

\end{remark}

\end{section}

\bigskip
\bigskip
\goodbreak

\begin{section}{The Hausdorff Series.}

Let $e^x=\exp x=1+\sum_{t\in M(x)} a(t)t\in K\{\{x\}\}$ be the
unique series with constant term 1 such that
$\exp(x)\exp(x)=\exp(2x)$, $\exp'(x)= \exp(x)$, see \cite{lidg}.

The composition inverse of $\exp x - 1$ is the series
$\log(1+x)=\sum_{t\in M(x)} b(t)t$.

We write $\exp(x)\exp(y) = 1 + \sum_{s\in M(x,y)} d(s)s$. The
Hausdorff series $H(x,y)= \sum_{\tau\in M(x,y)} c(\tau)\tau$ is
defined by $H(x,y)=\log(\exp(x)\exp(y))$.

Let $H=\sum_{n=1}^{\infty}H_n$, $H_n$ homogeneous of degree $n$.
Then $H_1=x+y$, $H_2=\frac{1}{2}(xy-yx)$.

\medskip

Similar to the classical case of associative variables, the
Hausdorff series has the following property.

\begin{thm}
The components $H_n$ of the non-associative Hausdorff series
$H(x,y)$ are primitive elements for the co-addition.
\end{thm}
\begin{proof}
\begin{itemize}
\item[1)]
 We first show, that $e^{ x\otimes 1+1\otimes x}=e^x\otimes e^x$.
\newline
Let $f(x)=e^{ x\otimes 1+1\otimes x}$. Substitution of $ x\otimes
1+1\otimes x$ for $x$ in $e^{2x}=e^xe^x$ yields $f(2x)= e^
{x\otimes 1+1\otimes x}e^{ x\otimes 1+1\otimes x}=f(x)f(x)$. Let
$g(x)=e^x\otimes e^x$. Clearly $g(x)g(x)=g(2x)$.
\newline
Writing $f=\sum f_n$ and $g=\sum g_n$ as sums over homogeneous
elements, one has $f_0=g_0=1\otimes 1, f_1=g_1=x\otimes 1+1\otimes
x$. Then $f_n$, and similarly $g_n$, are uniquely determined by
the equation $(2^n-2)f_n(x)= \sum_{0< i< n}f_i(x)\cdot
f_{n-i}(x)$, compare also \cite{lidg}. Thus $f=g$.
\newline
It follows that
$\Delta(e^xe^y)=\Delta(e^x)\Delta(e^y)=e^{\Delta(x)}e^{\Delta(y)}=e^xe^y\otimes
e^xe^y$.
\item[2)]
The rest of the proof goes along the well-known line, cf.\
\cite{lireu} \S 3.
\newline
To $\Delta( H(x,y) ) = H(x,y)\otimes 1+ 1\otimes H(x,y)$ we can
apply $\exp(x)-1$, to get the equivalent equation
\begin{equation*}
\exp(\Delta(\log(\exp(x)\exp(y))))-1=\exp(H\otimes 1+1\otimes
H)-1.
\end{equation*}
Writing $Z(x,y)$ for the series $\exp(x)\exp(y)-1$ without
constant term, the lefthandside is given by
$e^{\Delta(\log(1+Z(x,y)))}-1= \Delta\bigl(
e^{\log(1+Z(x,y))}-1\bigr)=\Delta(Z(x,y))=\Delta(e^xe^y)-1$. Since
$e^{H\otimes 1+1\otimes H}-1= e^xe^y\otimes
e^xe^y-1=\Delta(e^xe^y)-1$ by 1), the equation is true.
\end{itemize}
\end{proof}

\begin{defn}
Let $\vert\ \vert: M(x,y)\to M(z)$ be the homomorphism given by
$x\mapsto z, y\mapsto z$. We call $\vert\tau\vert$ the underlying
(unlabeled) tree of $\tau\in M(x,y)$.
 Let $t\in M(z)$ be of degree $n$, and
let $s_1,...,s_n\in M(x,y)$. Then $t(s_1,...,s_n)$ denotes the
result of grafting each $s_i$ to the $i$-th leaf of $t$.
\end{defn}

\begin{proposition}
The coefficients $c(\tau)$, $\tau \in M(x,y)$, of the Hausdorff
series are given by
\begin{equation*}
c(\tau)=d(\tau)-\sum_{k=2}^{\deg\tau}c_k(\tau)
\end{equation*}
\begin{equation*}
\text{with } c_k(\tau)=\begin{cases} c(\tau) &:\ k =1\\
a(\vert\tau\vert) &:\ k = \deg\tau
\\
\frac{1}{2^k-2}\sum_{l=1}^{k-1}c_l(\tau_1)c_{k-l}(\tau_2) &:\
2\leq k < \deg\tau, \text{ and } \tau=\tau_1\tau_2.
\end{cases}
\end{equation*}
 while $d(\tau)$ is the coefficient of $\tau$ in $\exp(x)\exp(y)$, given by
\begin{equation*}
 d(\tau)=\begin{cases}
a(\tau) &\text { if } \tau\in M(x) \text{ or } \tau\in M(y),\\
 a(\tau_1)a(\tau_2)  &\text{ if } \tau=\tau_1\tau_2, 1\neq \tau_1\in M(x),
 1\neq \tau_2\in M(y),\\
0 &\text{ else.}\\
\end{cases}
\end{equation*}
\end{proposition}
\begin{proof}
Insertion of $H(x,y)$ into $\sum_{t\in M(x)} a(t)t$ yields
\begin{equation*}
e^{H(x,y)}=1+\sum_{k=1}^{\infty}\sum_{t,s_1,...,s_k}
a(t)c(s_1)\cdot...\cdot c(s_n)t(s_1,...,s_k),
\end{equation*}
 where the second sum is
over all $t\in M(z)$ with $\deg t=k$ and $s_1,..,s_k\in M(x,y)$.
Now we can compare the coefficients of
\begin{equation*}
e^{H(x,y)}=1+\sum_{k=1}^{\infty}\sum_{t,s_1,...,s_k \choose
t(s_1,...,s_k)=\tau} a(t)c(s_1)\cdot...\cdot c(s_k)\tau,
\end{equation*}
with the coefficients of $\exp(x)\exp(y)$, which are easily
determined to be the values given above. Using $c(x)=a(x)=1$, we
get that
\begin{equation*}
c(\tau)=d(\tau)-a(\vert\tau\vert)-\sum_{k=2}^{\deg\tau-1}\sum_{\tau=t(s_1,...,s_k)}a(t)c(s_1)\cdot...\cdot
c(s_k)=d(\tau)-\sum_{k=2}^{\deg\tau}c_k(\tau).
\end{equation*}

\end{proof}

\begin{remark}
The $c(\tau)$ can be recursively computed as stated above, because
the coefficients $a(t)$ of $\exp(x)$ satisfy the formula
\begin{equation*}
a(t)=\frac{a(t_1)a(t_2)}{2^n-2}, \text{ if } t=t_1\cdot t_2,\ \
n:=\deg(t).
\end{equation*}
\end{remark}

\begin{example}

For $\tau\in M(x)$, $\deg\tau >1$, $c(\tau)=0$, as
$d(\tau)=a(\vert\tau\vert)$, and because $c(s_i)=0$ for at least
one $i$ in the formula above.

The coefficients $c(x\cdot t_2)$, for $t_2\in M(y)$ of degree
$n-1$, are given by $a(x)a(t_2)-a(\vert x
t_2\vert)=a(t_2)(1-\frac{1}{2^n-2})=\frac{2^n-3}{2^n-2}a(t_2).$

\end{example}

\begin{lemma}
\begin{itemize}
\item[(i)]
There is a unique continuous involution $*$ given on
$K\{\{x,y\}\}$ such that
\begin{equation*}
x^*=y,  \ y^*=x.
\end{equation*}
\item[(ii)] It holds that $(\exp(x))^*=\exp(y)$.
\item[(iii)] $H(x,y)^*=H(x,y)$. Furthermore $H(x,y)_n^*=H(x,y)_n$
for all $n$.
\end{itemize}
\end{lemma}
\begin{proof}
For the involution $*$, we want to show that
$(\exp(x))^*=\exp(y)$. Then also
$H(x,y)^*=\log(\exp(x)\exp(y))^*=\log(\exp(x)\exp(y))=H(x,y)$.

Since $(\exp(x))^*=g(y)$ is a series in $y$ (with constant term 1)
satisfying $g(y)g(y)=g(2y)$, we conclude that $g(y)=\exp(y)$.

\end{proof}

\begin{example}

For $X=\{x,y\}$, the vector space of elements of multi-degree
$(2,1)$ in $A_0$ has $(x,x,y), (x,y,x), (y,x,x),$ together with
$[x,[x,y]]$ as a basis, see Example (\ref{extwovar}). The elements
$(x,x,y)^*=-(x,y,y), (x,y,x)^*=-(y,x,y), (y,x,x)^*=-(y,y,x),$
together with $[x,[x,y]]^*=-[y,[x,y]]$ form a basis of the
multi-degree $(1,2)$-part.

\bigskip

Further computation shows:
\end{example}

\begin{proposition}
The homogeneous part $H_3$ of degree 3 is given by
\begin{equation*}
\begin{split}
H_3=&\frac{1}{12}([x,[x,y]]+[x,[x,y]]^*)+\frac{1}{12}((y,x,x)+(y,x,x)^*)\\
&+\frac{5}{12}((x,x,y)+(x,x,y)^*)+
\frac{1}{4}((x,y,x)+(x,y,x)^*)\\
=&\frac{1}{3}\Bigl((x,x,y)+(x,y,x)+\frac{1}{4}\bigl([x,[x,y]]+[x^2,y]-x[x,y]-[x,y]x\bigr)\Bigr)\\
&+\frac{1}{3}\Bigl((x,x,y)+(x,y,x)+\frac{1}{4}\bigl([x,[x,y]]+[x^2,y]-x[x,y]-[x,y]x\bigr)\Bigr)^*.\\
\end{split}
\end{equation*}
\end{proposition}

\end{section}

\begin{remark}
 In this article we have considered the free magma algebra; we
 dealt with one operation without relations.
  The free commutative magma algebra is also very interesting. For example, the similarly
defined Hausdorff series is given by
\begin{equation*}
H^a(x,y)=x+y+\frac{1}{3}\bigl((x,x,y)+(x,x,y)^*\bigr) \text{ plus
terms of order }\geq 4.
\end{equation*}
It seems also important to study non-associative Hausdorff series
in $m>2$ non-associative variables $x_1,...,x_m$ (cf.\ \cite{lilo}
for the associative case).

\end{remark}

\bigskip
\end{document}